\documentclass[reqno,11pt,a4paper]{amsart}
\usepackage{amsmath, amsfonts, amssymb, amsthm}
\usepackage{cases}\usepackage{extarrows}
\usepackage{float}
\usepackage{mathrsfs}
\usepackage[T1]{fontenc}
\usepackage{mathrsfs}
\usepackage{enumitem}
 \usepackage{indentfirst}

\newtheorem{theorem1}{Theorem}[section]

\newtheorem{definition1}{Definition}[section]

\numberwithin{equation}{section}
\usepackage[%
colorlinks=true,  % 启用彩色链接而非边框
linkcolor=blue,  % 内部链接颜色
citecolor=blue,  % 引用链接颜色
urlcolor=black,   % URL链接颜色
anchorcolor=black,% 锚点链接颜色
]{hyperref}

\begin{document}

\title{Global dynamics and asymptotic stability for a nonlocal LANS-$\alpha$-type system}

\author[Guan]{Chunxia Guan}
\address{Chunxia Guan \newline
Department of Mathematics, Guangdong University of Technology, Guangzhou 510520, China}
\email{guanchunxia123@163.com}

\author{Kai Yan\textsuperscript{*}}
\thanks{\noindent $^*$Corresponding author.}
\address{Kai Yan (Corresponding author)\newline
School of Mathematics and Statistics\\
Huazhong University of Science and Technology, Wuhan 430074,  China}
\email{kaiyan@hust.edu.cn}

\author[Zhang]{Qianyuan Zhang}
\address{Qianyuan Zhang\newline
	School of Mathematics and Statistics\\
	Huazhong University of Science and Technology, Wuhan 430074,  China}
\email{qianyuanzhang@hust.edu.cn}

\begin{abstract}
This paper is devoted to the global dynamics of a nonlocal Lagrangian-averaged Navier-Stokes-$\alpha$ (LANS-$\alpha$)-type system (also known as the viscous Euler-Poincar\'e system) in one, two, and three spatial dimensions. In the unforced periodic setting with zero mean, we establish global well-posedness for the initial-value problem and prove that the corresponding solution semigroup is dissipative and possesses a compact global attractor consisting of the singleton $\{0\}$. Furthermore, we obtain
the following global-in-time algebraic asymptotic stability estimate:
\[
\|u(t,\cdot)\|_{H^3}\lesssim t^{-1/2},\qquad \forall\, t>0,
\]
which is derived from higher-order energy inequalities combined with a uniform Gronwall-type argument. Our results provide a complete description of the long-time behavior of the unforced periodic dynamics for this nonlocal LANS-$\alpha$-type system.
\end{abstract}

\maketitle


\vskip 0.2cm


\setcounter{equation}{0}

\section{Introduction}
In this paper, we study the following nonlocal Lagrangian-averaged Navier-Stokes-$\alpha$ (LANS-$\alpha$)-type system (also known as viscous Euler-Poincar\'e system), namely,
 \begin{equation}\label{eq1}
	\left\{\begin{array}{ll}
		\partial_t m+ u\cdot\nabla m- \epsilon\Delta m
		+ (\nabla u)^{T}\cdot m+ (\text{div} u) m=0,\\
		m=(I-\alpha^2\Delta)u,
	\end{array}\right.
\end{equation}
where one denotes $(u\cdot\nabla m)_i=u_j\partial_j m_i$,
$\left((\nabla u)^{T}\cdot m\right)_i=m_j\partial_i u_j $, and we have used the Einstein summation convention over repeated indices which will be obeyed throughout the whole paper. The velocity $u=u(t,x)$ and momentum $m=m(t,x)$ are defined from $\mathbb{R^+}\times \mathbb{A}^d$ to $\mathbb{R}^d$,
and $\mathbb{A}^d= \mathbb{R}^d$ or $\mathbb{T}^d$ with the torus $\mathbb{T}^d \triangleq \mathbb{R}^d /\mathbb{Z}^d$,
$d\in \mathbb{N}_+$ is the spacial dimensions, and $\epsilon>0$ is the coefficient of viscosity. The real constant $\alpha $ is dispersion parameter. For the sake of simplicity, we may take $\alpha=1.$ The Helmholtz relation $m=(I-\alpha^2\Delta)u$
is essential to the nonlocal character of the system when it is expressed in terms of $u$.

When taking $\epsilon=0$, system \eqref{eq1} becomes the corresponding Euler-Poincar\'{e} equations
%proposed by Holm \emph{et al.}  in \cite{MR1627802,MR2552212}, namely,
\begin{equation}\label{1mch}
m_t+u\cdot\nabla m+(\nabla u)^{T}\cdot m+ (\text{div} u)m=0,\qquad
m=(I-\alpha^2\Delta)u.
\end{equation}
The Euler-Poincar\'{e} system \eqref{1mch} associated with the diffeomorphism group in \cite{MR2103008}
was firstly proposed exactly in the way that
a class of its singular solutions generalize the peakon solutions of the classical Camassa-Holm  equation
to higher spatial dimensions \cite{MR2031278}.
Let us consider the  $d$-dimensional Lagrangian $L(u)$ defined by
\begin{equation}\label{1-Lagrangian}
L(u) \triangleq
\frac{1}{2} \int_{\mathbb{A}^d} u \cdot (I-\alpha\Delta)u dx.
\end{equation}
It is shown in \cite{MR1627802} that the corresponding semidirect-product Euler-Poincar\'{e} equations can be written as
\begin{equation}\label{1-EPDiff}
\frac{\partial}{\partial t} \frac{\delta L}{\delta u}
= - \pounds_{u}  \frac{\delta L}{\delta u},
\end{equation}
where $\pounds_{u}  \frac{\delta L}{\delta u}$ is the corresponding Lie derivatives of
$\frac{\delta L}{\delta u}$.
Then, in view of the Hamiltonian principle
$\delta \int_{t_0}^{t_1} L(u) dt=0$,
one evaluates the variational derivatives $\frac{\delta L}{\delta u}$
of the Lagrangian  \eqref{1-Lagrangian}, and then substitutes them into the Euler-Poincar\'{e} equations (\ref{1-EPDiff}),
which eventually recovers the system \eqref{1mch}. On the other hand, the system \eqref{1mch} can be viewed as a framework for modeling and analyzing fluid dynamics,
particularly for nonlinear shallow water waves, geophysical fluids and turbulence modeling \cite{MR2552212}.
The last three terms in the first system in \eqref{1mch} model convection, stretching and  expansion
of a fluid with velocity $u$ and momentum $m$, respectively.
In addition, its applications in imaging science were also discussed in \cite{MR2552212}.
The local well-posedness,  global existence and blow-up of solutions of the Cauchy problem for system \eqref{1mch}
have been discussed in \cite{MR2964772,MR3116009,MR4042343,MR3277196}.

In the special case, when the dimension $ d = 1, $  system \eqref{1mch} becomes the following famous Camassa-Holm (CH) equation:
\begin{equation*}
m_{t}+u m_x+2u_x m =0,\qquad
m=u-\alpha^2 u_{xx},
\end{equation*}
which is proposed as a bi-Hamiltonian equation \cite{MR636470} and models the unidirectional propagation of shallow water waves over a flat bottom \cite{MR1234453}
(see also a rigorous justification in shallow water \cite {MR2481064}).
It has a bi-Hamiltonian structure and is completely integrable \cite{MR1234453}.
Its solitary waves are peaked solitons (peakons) \cite{MR1234453,CAMASSA19941,MR2318158} (see also \cite{MR2257390,MR2753609} for the waves of great height in irrotational water waves), and they are orbitally stable \cite{MR1854962,MR1737505}.
The Cauchy problem  for the CH equation have been studied extensively.
It has been shown that this equation is locally well-posed \cite{MR1631589,MR1604278,MR1827098}.
Moreover, it has both global strong solutions \cite{MR1775353,MR1631589,MR1604278} and blow-up solutions within finite time \cite{MR1775353,MR1631589,MR1668586,MR1604278}
which is featured as  wave breaking \cite{CAMASSA19941,MR1668586} (namely, the wave remains bounded while its slope becomes unbounded in finite time \cite{MR483954}).
In addition, it possess global weak solutions, see the discussions in \cite{MR2278406,MR4108918,MR1773414}, for instance.

For $\epsilon>0$, system \eqref{eq1} is a viscous Euler-Poincar\'e system closely related to LANS-$\alpha$ model (also known as the viscous Camassa-Holm equations) which provides a regularized model of the classical Navier-Stokes equations. The LANS-$\alpha$ equations introduce a length scale through a Helmholtz relation between momentum and velocity and have been studied from the viewpoints of global regularity, long-time dynamics, and fluid modeling. In particular, Foias, Holm and Titi established global-in-time regularity for the three-dimensional LANS-$\alpha$ system and obtained estimates for the dimensions of its global attractor; see \cite{MR1837927,MR1878243}. These works provide an important background for the present study, although our objective is different: we focus on the unforced zero-mean periodic dynamics and completely identify the global attractor together with a quantitative algebraic convergence rate. In addition, non-autonomous LANS-$\alpha$ models with oscillating external forces have been studied from the viewpoint of uniform global attractors; see \cite{MR2754280}. The present autonomous, unforced, zero-mean problem is structurally different, and our main dynamical conclusion is the collapse of the global attractor to the single equilibrium $\{0\}$ together with an explicit algebraic convergence rate. We also note that algebraic large-time decay for viscous Camassa-Holm equations has been studied in the whole space, where the decay mechanism is tied to the dispersive and viscous structures, as well as the spatial localization of the data \cite{MR2438775,MR2457817}. The present periodic zero-mean problem is different in that the Poincar\'e inequality supplies a spectral gap, while the higher-order nonlocal energy estimates lead to the quantitative decay stated in Theorem~\ref{thm:main} below.

The purpose of this paper is to determine the global asymptotic behavior of system \eqref{eq1} in the unforced zero-mean periodic setting. We establish global well-posedness for $d=1,2,3$, construct the associated dissipative semigroup, and prove the existence of a compact global attractor. More importantly, we completely identify the attractor and show that it consists of the singleton $\{0\}$. We complement this qualitative result with a quantitative global-in-time algebraic decay estimate in $H^3(\mathbb{T}^d)$, which is obtained from higher-order energy inequalities combined with a uniform Gronwall-type argument. Thus, our contribution is not merely the existence of a global attractor, but its complete identification together with a quantitative description of the long-time convergence.

Now, we state the main result of the present paper as follows.

\begin{theorem1}\label{thm:main}
	Let $u_0\in H^s(\mathbb{A}^d)$ with $s>1+\frac{d}{2}$ and $d=1,2,3$. Then 
	\begin{enumerate}[label=(\arabic*), ref=(\arabic*)]
		\item\label{thm:main:a} the Cauchy problem for system \eqref{eq1} has a unique global solution $u\in C \big([0,\infty);H^s(\mathbb{A}^d)\big)\cap C^1\big([0,\infty);L^2(\mathbb{A}^d)\big)$. Moreover, for all $t>0$, we have
		\begin{equation*}
			\|u(t,\cdot)\|^2_{H^1(\mathbb{A}^d)}+2\epsilon\int_0^t\int_{\mathbb{A}^d}(|\nabla u|^2+|\Delta u|^2)(s,x)dxds=\|u_0\|^2_{H^1(\mathbb{A}^d)}.
		\end{equation*}
		\item\label{thm:main:b} the solution operator $S(t)u_0 \triangleq u(t)$ defines a continuous dissipative semigroup on $X_3(\mathbb{T}^d)\triangleq\big\{u\in H^3(\mathbb{T}^d)|\int_{\mathbb{T}^d}u(x)dx=0\big\}$, and possesses a compact global attractor which consists of the singleton $\{0\}$. Furthermore, there exists a constant $C=C(d,\epsilon,\alpha,\|u_0\|_{H^3(\mathbb{T}^d)})>0$ such that 
		\begin{equation}\label{eq2}
			\|u(t)\|_{H^3(\mathbb{T}^d)}\leq Ct^{-\frac{1}{2}},\quad\forall\ t>0,
		\end{equation}
		%Consequently, 
	%	\[ \lim_{t\to\infty}\|u(t)\|_{H^3(\mathbb{T}^d)}=0, \]
		showing that the equilibrium $u=0$ is globally asymptotically stable in $X_3(\mathbb{T}^d)$. 
	\end{enumerate}
\end{theorem1}

%\begin{remark1}
%Theorem~\ref{thm:main} gives two complementary descriptions of the long-time dynamics. The single equilibrium $\{0\}$ completely characterizes the global attractor, while the global-in-time estimate \eqref{eq2} gives an explicit algebraic decay rate at which trajectories approach this unique asymptotic state, which is obtained from higher-order energy inequalities combined with a uniform Gr\"onwall argument. Thus the result goes beyond the mere existence of a global attractor. Moreover, the global-in-time algebraic estimate \eqref{eq2} is obtained from higher-order energy inequalities combined with a uniform Gr\"{o}nwall argument.
%\end{remark1}

The remainder of the paper is organized as follows. In Section 2, for all $d\in\mathbb{N}_+$, we establish the local well-posedness for system \eqref{eq1} via Kato's semigroup theory. In Section 3, we prove global well-posedness in both $\mathbb{R}^d$ and $\mathbb{T}^d$ with $d=1,2,3$. In Section 4, we derive the continuous dissipative semigroup and show the complete characterization of the global attractor, as well as the global-in-time algebraic asymptotic stability estimate.

\section{Local Well-posedness}
\newtheorem {remark2}{Remark}[section]
\newtheorem{theorem2}{Theorem}[section]
\newtheorem{lemma2}{Lemma}[section]
\newtheorem{proposition2}{Proposition}[section]
\newtheorem{corollary2}{Corollary}[section]

In this section, we  establish the local well-posedness for the Cauchy problem of \eqref{eq1} in $\mathbb{A}^d$ with $d\in \mathbb{N}_+$.  
%the following viscous Eular-Poincar\'{e} equations:
% \begin{equation}\label{eq1}
%\left\{\begin{array}{ll}
%\partial_t m+ u\cdot\nabla m+ (\text{div} u) m
%+ (\nabla u)^{T}\cdot m=\epsilon\Delta m,\\
%m=(I-\Delta)u.
%\end{array}\right.
%\end{equation}
For this, let us first  set up the nonlocal transport equations form  for system \eqref{eq1}.
Indeed, we deduce from the system \eqref{eq1} that
\begin{align}\label{1-6-1}
  (I- \Delta) (\partial_t u+u\cdot \nabla u)
&= \partial_t m+u\cdot \nabla u- \Delta (u\cdot \nabla u) \nonumber\\ 
&= u\cdot \nabla (\Delta u)-\Delta(u\cdot \nabla u)+ \nabla u^T\cdot(\Delta u)+(\Delta u)\text{div}u \nonumber\\  
&\quad -u\nabla u-u(\text{div}u)+\epsilon(I- \Delta)\Delta u.
\end{align}

On the other hand,  the following three identities hold true:
\begin{equation}\label{1-6-3}
\quad\quad  u\cdot \nabla (\Delta u)-\Delta(u\cdot \nabla u)
=-\text{div}(\nabla u\nabla u+\nabla u\nabla u^T)+(\nabla u)\cdot \nabla(\text{div}u),
\end{equation}
\begin{equation}\label{1-6-4}
\nabla u^T\cdot(\Delta u)
=\text{div}(\nabla u^T \nabla u-\frac{1}{2}|\nabla u|^2 I ),
\end{equation}
\begin{equation}\label{1-6-5}
(\nabla u)\cdot \nabla(\text{div}u)+(\Delta u)\text{div}u
=\text{div}(\nabla u(\text{div}u)).
\end{equation}

Hence, in view of (\ref{1-6-1})-(\ref{1-6-5}), the  initial-value problem for system \eqref{eq1} can be rewritten as follows:
\begin{equation}\label{1mchnew}
\left\{\begin{array}{ll}
\partial_t u+u\cdot \nabla u= F(u)+\epsilon\Delta u, &(t, x)\in \mathbb{R}^+\times \mathbb{A}^d,\\
u(0,x)= u_{0}(x), &x\in \mathbb{A}^d,
\end{array}\right.
\end{equation}
where
\begin{align}\label{1-F1}
 F(u)\triangleq &-(I-\Delta)^{-1}\text{div}(\nabla u\nabla u+\nabla u\nabla u^T -\nabla u^T \nabla u-\nabla u (\text{div} u)
+\frac{1}{2}|\nabla u|^2 I)\nonumber \\
&-(I-\Delta)^{-1} (u(\text{div} u)+ u\nabla u).
\end{align}

Next, we will apply Kato's semigroup theory to establish the local
well-posedness for the Cauchy problem \eqref{1mchnew}.
For convenience, we state here Kato's theory in the form suitable
for our purpose. Consider the abstract quasi-linear evolution equation:
\begin{equation}\label{eq3}
\frac{dv}{dt}+A(v)v=f(v),\ \ t>0,\ \ v(0)=v_{0}.
\end{equation}

Let $X$ and $Y$ be Hilbert spaces such that $Y$ is continuously and
densely embedded in $X$ and let $Q : Y \rightarrow X$ be a
topological isomorphism. Let $L(Y,X)$ denote the space of all bounded
linear operators from $Y$ to $X$
($L(X)$, if $X=Y$.). Assume that \\
($i$) $A(y)\in L(Y,X)$ for $y\in Y$ with
$$ \left\|(A(y)-A(z))w\right\|_{X}\leq\mu_{1}\left\|y-z\right\|_{X}\left\|w\right\|_{Y},\ \ \ \ y,z,\omega\in Y,$$
and $A(y)\in G(X,1,\beta)$, (i.e. $A(y)$ is quasi-$m$-accretive),
uniformly on bounded sets in $Y$.\\
($ii$) $QA(y)Q^{-1}=A(y)+B(y)$, where $B(y)\in L(X)$ is bounded,
uniformly on bounded sets in $Y$. Moreover,
$$ \left\|(B(y)-B(z))w\right\|_{X}\leq\mu_{2}\left\|y-z\right\|_{Y}\left\|w\right\|_{X},\ \ \ \ y,z\in Y ,
\ \omega\in X.$$ 
($iii$) $f: Y\rightarrow Y$ and extends also to a map from $X$ to $X$.
$f$ is bounded on bounded sets in $Y$, and
$$ \left\|f(y)-f(z)\right\|_{Y}\leq \mu_{3}\left\|y-z\right\|_{Y}, \ \ \ y,z\in Y,$$
$$
\left\|f(y)-f(z)\right\|_{X}\leq\mu_{4}\left\|y-z\right\|_{X}, \ \ 
y,z\in Y.$$ Here $\mu_{1}, \mu_{2}, \mu_{3}$ and $ \mu_{4}$ depend
only on max$\{\left\|y\right\|_{Y},
\left\|z\right\|_{Y} \}$.

 \begin{lemma2}\label{T2-1}
 \cite{MR407477} Assume that (i), (ii) and (iii) hold. Given  $ v_{0}\in Y$,  there
is a maximal $T>0$ depending only on $\|v_{0}\|_{Y}$
and a unique solution v to system \eqref{eq3} such that
$$
v=v(\cdot,v_{0})\in C([0,T);Y)\cap C^{1}([0,T);X).
$$
Moreover, the map $v_{0}\mapsto v(\cdot,v_{0})$ is continuous
from Y to
$
 C([0,T);Y)\cap C^{1}([0,T);X).
$
 \end{lemma2}

Before proving the main result of this section, we recall some useful lemmas as follows.

\begin{lemma2}\label{L2-1}
\cite{MR407477} Let $r,t$ be real numbers such that $-r<t\leq r.$ Then
\begin{equation*}
	\|fg \|_{H^{t}(\mathbb{A}^d)}\leq
	c\|f\|_{H^{r}(\mathbb{A}^d)}\|g\|_{H^{t}(\mathbb{A}^d)},  \ \ \ \ if \
	r>\frac{d}{2},
\end{equation*}
\begin{equation*}
	\|fg \|_{H^{t+r-\frac{d}{2}}(\mathbb{A}^d)}\leq
	c\|f\|_{H^{r}(\mathbb{A}^d)}\|g\|_{H^{t}(\mathbb{A}^d)}, \ if \
	r<\frac{d}{2},
\end{equation*}
where c is a positive  constant depending on r, t and d.
\end{lemma2}

\begin{lemma2}\label{L2-2}
 \cite{MR407477} Let $f\in H^s(\mathbb{A}^d), s>\frac{d}{2}+1.$ Then
 \begin{equation*}
 	\|[\Lambda^s,
 	M_f]\Lambda^{1-s}\|_{L(L^2(\mathbb{A}^d))}\leq
 	c\|\nabla f\|_{H^{s-1}(\mathbb{A}^d)},
 \end{equation*}
where $\Lambda=(I-\Delta)^{\frac{1}{2}}$, 
$M_f$ is the operator of multiplication by $f$ and $c$ is a
positive constant, $\|\cdot\|_{L(L^2(\mathbb A^d))}$ denotes the operator norm in $L^2(\mathbb A^d)$.
\end{lemma2}

The main result of this section is the following theorem.

\begin{theorem2}\label{T2-2}
Given $u_0\in H^{s}(\mathbb{A}^d),$ $s>1+\frac{d}{2}$ $(d\in\mathbb N_+)$, there
exists a maximal $T=T(\|u_{0}\|_{H^{s}(\mathbb{A}^d)})>0$, and a unique solution $u$ to
the system \eqref{1mchnew} such that
\begin{equation*}
	u=u(\cdot,u_{0})\in C\big([0,T); H^{s}(\mathbb{A}^d)\big)\cap C^1\big([0,T); L^{2}(\mathbb{A}^d)\big).
\end{equation*}
Moreover, the map $u_0\mapsto u(\cdot,u_0)$ is continuous from $H^s(\mathbb{A}^d)$ to $C\big([0,T); H^{s}(\mathbb{A}^d)\big)\cap C^1\big([0,T); L^{2}(\mathbb{A}^d)\big)$. 
\end{theorem2}
\begin{proof}
Let $A(u)\triangleq u\cdot\nabla-\epsilon\Delta $
and $f(u)\triangleq F(u)$ given by \eqref{1-F1}.
Set $Y=H^{s}(\mathbb{A}^d),\
X=L^2(\mathbb{A}^d), \
\Lambda=(1-\Delta)^{\frac{1}{2}}$ and $Q=\Lambda^s $. Obviously, Q is an isomorphism of $H^{s}(\mathbb{A}^d)$ onto $L^2(\mathbb{A}^d)$. In order to prove Theorem \ref{T2-2}. In view of Lemma \ref{T2-1}, we only need to verify that $A(u)$ and $f(u)$ satisfy the conditions $(i)-(iii)$.

 Similar as the arguments in \cite{MR2823888}, we have $A(u)\in G(X, 1, \beta(u)).$

Now we prove that $A(y)\in L(Y,X)$ for
$y\in Y$, and 
$$ \|(A(y)-A(z))w\|_{X}\leq\mu_{1}\|
y-z\|_{X}\|w\|_{Y},$$ for $y,z,w\in Y.$
Indeed, by the definition of $A(u)$ and the fact that $H^t(\mathbb A^d)\hookrightarrow L^\infty(\mathbb A^d)$ as $t>\frac{d}{2}$, we deduce that for all
$z,w\in Y$,
\begin{align*}
 \|A(z)w\|_{X}\leq\|z\cdot\nabla w\|_{X}+\epsilon\|\Delta w\|_{X}\leq\|z\|_{L^\infty(\mathbb A^d)}\|\nabla w\|_{X}+\epsilon\|\Delta w\|_{X}\leq(\|z\|_{Y}+\epsilon)\| w\|_{Y}.
\end{align*}
 On the other hand, for all $y,z,w\in Y,$
\begin{align*}
 \|(A(y)-A(z))w\|_{X} \leq\|\nabla w\|_{L^{\infty}(\mathbb{A}^d)}\|y-z\|_{X} \leq\|y-z\|_{X}\|w\|_{Y}.
\end{align*}

Note that $B(z)w\triangleq Q A(z)Q^{-1}w-A(z)w=[\Lambda^s, z\cdot\nabla]\Lambda^{-s} w$.
Thus, for $w\in X$ and $y,z\in
Y$,  we infer that
\begin{align*}
\big\|(B(y)w-B(z)w)\big\|_{X}
=&\big\|[\Lambda^s,(y-z)\cdot\nabla]\Lambda^{-s} w\big\|_{L^2}\\
\leq&\big\|[\Lambda^s,
y-z]\Lambda^{1-s} w\|_{L(L^2)}\|\Lambda^{-1}\nabla w\big\|_{L^2}\\
\leq& C\big\|\nabla(y-z)\big\|_{H^{s-1}}\big\|w\big\|_{L^2}\\
\leq& C\|y-z\|_{Y}\|w\|_{X},
\end{align*}
where we applied Lemma \ref{L2-2} in view of $s>1+\frac{d}{2}$.
Taking $z=0$ in the above inequality, we obtain $B(u)\in
L(X)$.

 Denote $g(u)\triangleq\nabla u\nabla u+\nabla u\nabla u^T -\nabla u^T \nabla u-\nabla u (\text{div} u)
+\frac{1}{2}|\nabla u|^2 I$. Then for every $u,v\in Y,$ we  have 
\begin{align*}
 &\|f(u)-f(v)\|_{Y}\\
 =&\|\text{div}(g(u)-g(v))+u(\text{div} u)+ u\nabla u-v(\text{div} v)- v\nabla v\|_{H^{s-2}(\mathbb{A}^d)}\\
 \leq&\|g(u)-g(v)\|_{H^{s-1}(\mathbb{A}^d)}+\|u(\text{div} u)+u\nabla u-v(\text{div} v)-v\nabla v \|_{H^{s-2}(\mathbb{A}^d)}\\
 \leq&\|\nabla u\nabla (u-v)+\nabla (u-v)(\nabla v)+\nabla u\nabla (u-v)^T+\nabla (u-v)\nabla v^T-\nabla (u-v)^T\nabla u\\
&-\nabla v^T\nabla (u-v)-\nabla u(\text{div}(u-v))-\nabla (u-v)(\text{div} v)+\frac{1}{2}(\nabla u:\nabla(u-v)+\nabla (u-v):\nabla v)I\|_{H^{s-1}(\mathbb A^d)}\\
&+\|u\text{div} (u-v)+(u-v)\text{div} v+u\nabla(u-v)+(u-v)\nabla v\|_{H^{s-1}(\mathbb{A}^d)}\\
\leq&C(\|\nabla u\|_{H^{s-1}(\mathbb{A}^d)}+\|\nabla v\|_{H^{s-1}(\mathbb{A}^d)}+\|u\|_{H^{s-1}(\mathbb{A}^d)}+\|\text{div} v\|_{H^{s-1}(\mathbb A^d)})\\
&(\|\nabla(u-v)\|_{H^{s-1}(\mathbb{A}^d)}+\|\text{div}(u-v)\|_{H^{s-1}(\mathbb{A}^d)}+\|u-v\|_{H^{s-1}(\mathbb{A}^d)})\\
\leq& C \|u-v\|_Y,
\end{align*}
where $C$ is a constant depending on $\|u\|_Y$ and $\|v\|_{Y}$, and we have used the fact that  $H^{s-1}(\mathbb A^d)$ is a Banach algebra with $s>1+\frac{d}{2}$.

Finally, we only need to verify that
\begin{align*}
 &\|f(y)-f(z)\|_X\\
 =&\|\text{div}(g(u)-g(v))+u(\text{div} u)+u\nabla u-v(\text{div} v)-v\nabla v\|_{H^{-2}(\mathbb{A}^d)} \\
 \leq&\|g(u)-g(v)\|_{H^{-1}(\mathbb{A}^d)}+\|u(\text{div}u) +u\nabla u-v\text{div} v-v\nabla v\|_{H^{-2}(\mathbb{A}^d)}\\
 \leq&\|\nabla u\nabla (u-v)+\nabla (u-v)\nabla v+\nabla u\nabla (u-v)^T+\nabla (u-v)\nabla v^T-\nabla (u-v)^T\nabla u-\nabla v^T\nabla (u-v)\\
&-\nabla u(\text{div}(u-v))-\nabla (u-v)(\text{div}v) +\frac{1}{2}(\nabla u:\nabla(u-v)+\nabla (u-v):\nabla v)I\|_{H^{-1}(\mathbb{A}^d)}\\
&+\|u(\text{div} (u-v))+(u-v)(\text{div}v) +u\nabla(u-v)+(u-v)\nabla v\|_{H^{-1}(\mathbb{A}^d)}\\
\leq&C \left(\|\nabla u\|_{H^{s-1}(\mathbb{A}^d)}+\|\nabla v\|_{H^{s-1}(\mathbb{A}^d)}+\|\text{div}v\|_{H^{s-1}(\mathbb{A}^d)}+\|u\|_{H^{s-1}(\mathbb{A}^d)}\right)\\
&(\|\nabla (u-v)\|_{H^{-1}(\mathbb{A}^d)}+\|\text{div}(u-v)\|_{H^{-1}(\mathbb{A}^d)}+\|u-v\|_{H^{-1}(\mathbb A^d)})\\
\leq &C\|u-v\|_X,
\end{align*}%
where we applied Lemma \ref{L2-1} with $r=s-1>\frac{d}{2}, t=-1$, and $C$ depending on $\|u\|_Y$ and $\|v\|_{Y}$. Therefore, we have completed the proof Theorem \ref{T2-2}.
\end{proof}

\section{Global well-posedness}
\newtheorem {remark3}{Remark}[section]
\newtheorem{theorem3}{Theorem}[section]
\newtheorem{lemma3}{Lemma}[section]
\newtheorem{proposition3}{Proposition}[section]
\newtheorem{corollary3}{Corollary}[section]

Our  attention in the present section will be focused on  the  existence of global strong solutions  to system \eqref{1mchnew} by proving Theorem \ref{thm:main} \ref{thm:main:a}.
Firstly, we give an energy estimate  as follows.

\begin{lemma3}\label{L3-1}
Let $u_{0}\in H^s(\mathbb{A}^d), s>1+\frac{d}{2}$ and assume that  $T>0$ is the
maximal existence time of the corresponding solution $u$
to system \eqref{1mchnew}. Then for all $t\in[0,T)$, we have
\begin{equation}\label{eq4}
	\|u(t,\cdot)\|^2_{H^1(\mathbb{A}^d)}+2\epsilon\int_0^t\int_{\mathbb A^d}(|\nabla u|^2+|\Delta u|^2)(s,x)dxds=\|u_0\|^2_{H^1(\mathbb{A}^d)}.
\end{equation}
\end{lemma3}

\begin{proof} 
	By a simple density argument, we only need to prove this lemma for some $s>3+\frac{d}{2}$.
Multiplying  by $2u$ on both sides of system \eqref{eq1} and integrating over
$\mathbb{A}^d$,  applying integration by parts,  we deduce
\begin{align*}
&\frac{d}{dt}\int_{\mathbb A^d}u\cdot m(t,x)dx=2\int_{\mathbb A^d}u\cdot m_t(t,x)dx\\
=&-2\int_{\mathbb A^d}u\cdot(u\cdot \nabla m)(t,x)dx-2\int_{\mathbb A^d}(\text{div}u)(u\cdot m)(t,x)dx\\
&-2\int_{\mathbb A^d}u\cdot(m(\nabla u))(t,x)dx+2\epsilon\int_{\mathbb A^d}u\cdot\Delta m(t,x)dx\\
=&-2\int_{\mathbb A^d}(u_ju_i\partial_im_j)dx-2\int_{\mathbb A^d}(\text{div}u)(u\cdot m)(t,x)dx\\
&-2\int_{\mathbb A^d}u\cdot(m(\nabla u))(t,x)dx+2\epsilon\int_{\mathbb A^d}\Delta u\cdot m(t,x)dx\\
=&2\int_{\mathbb A^d}(u_jm_j\partial_iu_i+u_im_j\partial_iu_j)(t,x)dx-2\int_{\mathbb A^d}(\text{div}u)(u\cdot m)(t,x)dx\\
&-2\int_{\mathbb A^d}u\cdot(m(\nabla u))(t,x)dx-2\epsilon\int_{\mathbb A^d}(|\nabla u|^2+|\Delta u|^2)(t,x)dx\\
=&-2\epsilon\int_{\mathbb A^d}(|\nabla u|^2+|\Delta u|^2)(t,x)dx.
\end{align*}
Integrating the above equality over $(0,t)$ with respect to the time variable, we obtain \eqref{eq4}.
\end{proof}

Next, we state the blow-up scenario for system \eqref{1mchnew}. To this end, let us recall the following useful lemma.

\begin{lemma3}\label{L3-2}
\cite{MR951744} If $r>0$, then there holds
$$\|[\Lambda^{r},f]g\|_{L^{2}(\mathbb{A}^d)}\leq C(\| \nabla f\|_{L^{\infty}(\mathbb{A}^d)}
\|\Lambda^{r-1}g\parallel_{L^2(\mathbb{A}^d)}+\|
\Lambda^r f\|_{L^2(\mathbb{A}^d)}\| g\|_{L^{\infty}(\mathbb{A}^d)}),
$$
where $C$ is a constant depending only on $r$ and $d$.
\end{lemma3}

\begin{lemma3}\label{L3-3}
Let $u_{0}\in H^{s}(\mathbb{A}^d)$, $s>1+\frac{d}{2}$ and
assume that $T>0$ is the maximal existence time of the corresponding
solution $u$ to system \eqref{1mchnew}. Then the
solution $u$
blows up $(\ i.e.,\ T<+\infty)$ if and only if
\[ \limsup_{t\rightarrow T}
(\|u(t,\cdot)\|_{L^\infty(\mathbb{A}^d)}+\|
\nabla u(t,\cdot)\|_{L^{\infty}(\mathbb{A}^d)})
=\infty. \]
\end{lemma3}

\begin{proof}
Applying the operator $\Lambda^{s}$ to the first equation in
system \eqref{1mchnew}, multiplying by $2\Lambda^{s} u$, and integrating over
$\mathbb{A}^d$, one gets
\begin{equation}\label{3.2}
\frac{d}{dt}\|
u\|^{2}_{H^{s}(\mathbb{A}^d)}=-2(u\cdot\nabla u,
u)_{s}-2(F(u),u)_{s}+2\epsilon(\Delta u,u)_{s},
\end{equation}
where
$(u,v)_s\triangleq\int_{\mathbb A^d}\Lambda^su\cdot\Lambda^svdx,$ $s\geq0.$
Thanks to Lemma \ref{L3-2} and integration by parts, we get
\begin{align}
|(u\cdot\nabla u,u)_{s}|&= |(\Lambda^s(u\cdot\nabla u),\Lambda^s u)_{0}|\nonumber\\
&= |([\Lambda^s,u\cdot]\nabla u,\Lambda^s
u)_{0}+(u\cdot\Lambda^s\nabla u,\Lambda^s u)_{0}|\nonumber\\
&\leq \|[\Lambda^s,u\cdot]\nabla u\|_{L^2(\mathbb{A}^d)}\|\Lambda^s
u\|_{L^2(\mathbb{A}^d)}+\frac{1}{2}|((\text{div}u)\Lambda^s u,\Lambda^s u)_{0}|\nonumber\\
& \leq (C\| \nabla u\|_{L^{\infty}(\mathbb{A}^d)}+\frac{1}{2}\|
\text{div}u\|_{L^{\infty}(\mathbb{A}^d)})\|u\|^2_{H^s(\mathbb{A}^d)}\nonumber\\
&\leq C\|\nabla u\|_{L^{\infty}(\mathbb{A}^d)}\|u\|^2_{H^s(\mathbb{A}^d)}.
\end{align}
While in view of the Cauchy-Schwarz inequality, one has
\begin{align}
|(F(u),u)_{s}|\leq\|F(u)\|_{H^s(\mathbb{A}^d)}\|u\|_{H^s(\mathbb{A}^d)}.
\end{align}
 Note that
 $-(I-\Delta)^{-1}\text{div},\, -(I-\Delta)^{-1} \in Op(S^{-1})$.
 Due to the Morse-type inequality, we infer 
 \begin{equation} \label{3.5}
 	\|F(u)\|_{H^s(\mathbb{A}^d)}\leq C(\|u\|_{L^\infty(\mathbb{A}^d)}+\|\nabla u\|_{L^\infty(\mathbb{A}^d)})\|u\|_{H^s(\mathbb{A}^d)}.
 \end{equation}
Since 
$2\epsilon(\Delta u,u)_{s}=-2\epsilon\|\nabla u\|^2_{H^s}\leq0$, it follows from (\ref{3.2})-(\ref{3.5}) that
\begin{equation}\label{eq15}
\frac{d}{dt}\|u\|^2_{H^s(\mathbb{A}^d)}\leq C(\|u\|_{L^\infty(\mathbb{A}^d)}+\|\nabla u\|_{L^\infty(\mathbb{A}^d)})\|
u\|_{H^s(\mathbb{A}^d)}^2.
\end{equation}
If there exists $M>0$ such that $\limsup\limits_{t\rightarrow T}
(\|u\|_{L^\infty(\mathbb{A}^d)}+\|\nabla u\|_{L^\infty(\mathbb{A}^d)})
\leq M,$  then \eqref{eq15} together with  Gronwall's inequality implies that the solution $u$ will not blow up within the lifespan.                          

On the other hand, by Sobolev's embedding theorem, we see that if
\[ \limsup_{t\rightarrow T}
(\|u(t,\cdot)\|_{L^{\infty}(\mathbb{A}^d)}+\|
\nabla u(t,\cdot)\|_{L^{\infty}(\mathbb{A}^d)})
=\infty, \]
then the solution will blow up in finite time. Therefore, we complete the proof of the lemma.
\end{proof}

\begin{lemma3}\label{L3-4}
Let $u_{0}\in H^s(\mathbb{A}^d)$ with $s>1+\frac{d}{2}$ and $1\leq d\leq3$. Assume that $T>0$ is the maximal existence time of the corresponding
solution $u$ to system (2.1). Then there exists some $C=C(\epsilon,\|u_0\|_{H^s(\mathbb A^d)})>0$ such that
\begin{equation}\label{eq2-5}
\|u(t,\cdot)\|_{L^{\infty}(\mathbb{A}^d)}+\|\nabla u(t,\cdot)\|_{L^{\infty}(\mathbb{A}^d)}\leq Ce^{Ce^{Ct}},\ \ \
\ \forall t\in[0,T).
\end{equation}
\end{lemma3}

\begin{proof}
 As before, we may assume $s>3+\frac{d}{2}$ to prove this lemma.  Multiplying by $m$ on both sides of system \eqref{eq1} and integrating over 
$\mathbb{A}^d$, using integration by parts and the H\"{o}lder inequality, one infers
\begin{align}\label{3-4}
&\quad\frac{1}{2}\frac{d}{dt}\int_{\mathbb A^d}|m|^2dx+\epsilon\int_{\mathbb A^d}|\nabla m|^2dx\nonumber\\
&=-\int_{\mathbb A^d}m\cdot(u\cdot \nabla m)dx-\int_{\mathbb A^d}(\text{div}u)|m|^2dx-\int_{\mathbb A^d}m\cdot((\nabla u)^T\cdot m)dx\nonumber\\
&=-\int_{\mathbb A^d}m\cdot(u\cdot \nabla m)dx+2\int_{\mathbb A^d}u_im_j\partial_im_jdx+\int_{\mathbb A^d}(m_iu_j\partial_im_j+u_jm_j\partial_im_i)dx\nonumber\\
&=\int_{\mathbb A^d}m\cdot(u\cdot \nabla m)dx+\int_{\mathbb A^d}m\cdot( (\nabla m)^T\cdot u)dx+\int_{\mathbb A^d}(\text{div}m)u\cdot mdx\nonumber\\
&\leq3\|u\|_{L^4(\mathbb{A}^d)}\|m\|_{L^4(\mathbb{A}^d)}\|\nabla m\|_{L^2(\mathbb{A}^d)}.
\end{align}
By Lemma \ref{L3-1}, for all $t\in [0,T)$, we have  $\|u(t,\cdot)\|_{H^1}\leq \|u_0\|_{H^1}.$ Then the Sobolev embedding theorem yields
\begin{equation}\label{eq5}
	\|u\|_{L^4(\mathbb{A}^d)}\leq C\|u\|_{H^1(\mathbb{A}^d)}\leq C\|u_0\|_{H^1(\mathbb{A}^d)}\ \ \text{for}\ 1\leq d\leq4.
\end{equation}
Thanks to the Gagliardo-Nirenberg inequality,  there exists some $\theta=1-\frac{d}{4}\in[0,1)$ such that
\begin{equation}
	\|m\|_{L^4(\mathbb{A}^d)}\leq C\|m\|^\theta_{L^2(\mathbb{A}^d)}\|\nabla m\|^{1-\theta}_{L^2(\mathbb{A}^d)}.
\end{equation}
 Then, in view of the following Young inequality 
 \begin{equation}\label{eq6}
 	ab\leq C(\epsilon)a^p+\epsilon b^q\ \ \text{with}\ \frac{1}{p}+\frac{1}{q}=1\ \ \text{and}\ \epsilon>0,
 \end{equation}
 one has
 \[ 3\|u\|_{L^4(\mathbb{A}^d)}\|m\|_{L^4(\mathbb{A}^d)}\leq C \|m\|_{L^2(\mathbb{A}^d)} +\frac{\epsilon}{2}\|\nabla m\|_{L^2(\mathbb{A}^d)}, \]
where $C=C(\epsilon,\|u_0\|_{H^1(\mathbb{A}^d)})>0$.  Thus, we get
\begin{align}
\frac{1}{2}\frac{d}{dt}\int_{\mathbb A^d}|m|^2dx\nonumber&\leq\nonumber C\|m\|_{L^2(\mathbb{A}^d)}\|\nabla m\|_{L^2(\mathbb{A}^d)} +\frac{\epsilon}{2}\|\nabla m\|^2_{L^2(\mathbb{A}^d)}-\epsilon\int_{\mathbb A^d}|\nabla m|^2dx\\
&\leq C\nonumber\|m\|^2_{L^2(\mathbb{A}^d)}+\frac{\epsilon}{2}\|\nabla m\|^2_{L^2(\mathbb{A}^d)} +\frac{\epsilon}{2}\|\nabla m\|^2_{L^2(\mathbb{A}^d)}-\epsilon\int_{\mathbb A^d}|\nabla m|^2dx\\
&\leq C\|m\|^2_{L^2(\mathbb{A}^d)},
\end{align}
which along with  Gronwall's inequality implies
\begin{equation}\label{3-6}
\|u(t,\cdot)\|^2_{H^2(\mathbb{A}^d)}=\|m\|^2_{L^2(\mathbb{A}^d)}\leq
e^{Ct}\|m_{0}\|^2_{L^2(\mathbb{A}^d)}\leq Ce^{Ct},
\end{equation}
where $C=C(\epsilon, \|u_0\|_{H^s(\mathbb A^d)})>0.$

Next, we estimate the norm $\|u\|_{H^3(\mathbb{A}^d)}$. Differentiating system \eqref{eq1} with respect to $x_k$,  multiplying by $\partial_km$ and integrating over
$\mathbb{A}^d$, then integration by parts and thanks to the H\"{o}lder inequality,  we deduce
\begin{align}\label{3-7}
&\frac{1}{2}\frac{d}{dt}\int_{\mathbb{A}^d}|\nabla m|^2dx+\epsilon\|D^2m\|_{L^2(\mathbb{A}^d)}\nonumber\\
=&\frac{1}{2}\frac{d}{dt}\int_{\mathbb A^d}|\nabla m|^2dx+\epsilon\int_{\mathbb A^d}\partial_j\partial_k m_i\cdot\partial_j\partial_k m_idx\nonumber\\
=&-\int_{\mathbb A^d}\partial_km_i\partial_k(u_j\partial_jm_i)dx-\int_{\mathbb A^d}\partial_km_i\partial_k((\text{div}u)m_i)dx-\int_{\mathbb A^d}\partial_km_i\partial_k(m_j\partial_iu_j)dx\nonumber\\
=&\int_{\mathbb A^d}\Delta m_i(u_j\partial_jm_i)dx+\int_{\mathbb A^d}\Delta m_i((\text{div}u)m_i)dx+\int_{\mathbb A^d}\Delta m_i(m_j\partial_iu_j)dx\nonumber\\
\leq& C\|\Delta m\|_{L^2(\mathbb A^d)}\left(\| u\|_{L^4(\mathbb{A}^d)}\|\nabla m\|_{L^4(\mathbb{A}^d)}+\|m\|_{L^4(\mathbb{A}^d)}\|\nabla u\|_{L^4(\mathbb{A}^d)}\right),
\end{align}
which together with \eqref{eq5}-\eqref{eq6} and \eqref{3-7} yields
\begin{align*}
&\quad\frac{1}{2}\frac{d}{dt}\int_{\mathbb A^d}|\nabla m|^2dx+\epsilon\|D^2m\|^2_{L^2(\mathbb A^d)}\\
&\leq\nonumber C\|\Delta m\|_{L^2(\mathbb A^d)}\left(\|u_0\|_{H^1(\mathbb A^d)}(C(\epsilon)\|\nabla m\|_{L^2(\mathbb A^d)}+\frac{\epsilon}{C}\|D^2m\|_{L^2(\mathbb A^d)})+\|\nabla u\|_{H^1(\mathbb A^d)}\|m\|_{H^1(\mathbb A^d)}\right)\\
&\leq\nonumber C\|D^2m\|_{L^2(\mathbb{A}^d)}\|\nabla m\|_{L^2(\mathbb{A}^d)}+\frac{\epsilon}{4}\|D^2m\|^2_{L^2(\mathbb{A}^d)}+C\|D^2m\|_{L^2(\mathbb{A}^d)}\|\nabla m\|_{L^2(\mathbb{A}^d)}\|m\|_{L^2(\mathbb{A}^d)}\\
&\leq\nonumber Ce^{Ct}\|\nabla m\|_{L^2(\mathbb{A}^d)}\|D^2m\|_{L^2(\mathbb{A}^d)}+\frac{\epsilon}{4}\|D^2m\|^2_{L^2(\mathbb{A}^d)}\\
&\leq\nonumber Ce^{Ct}\|\nabla m\|^2_{L^2(\mathbb{A}^d)}+\frac{\epsilon}{2}\|D^2m\|^2_{L^2(\mathbb{A}^d)}.
\end{align*}
This turns out that
\begin{equation*}
	\frac{d}{dt}\|\nabla m\|^2_{L^2(\mathbb{A}^d)}\leq Ce^{Ct}\|\nabla m\|^2_{L^2(\mathbb{A}^d)},
\end{equation*}
which along with Gronwall's inequality yields that 

\begin{equation}\label{eq7}
\|\nabla m\|^2_{L^2(\mathbb{A}^d)}\leq\|\nabla m_{0}\|^2_{L^2(\mathbb{A}^d)}
e^{C\int^t_0e^{C\tau} d\tau}\leq Ce^{Ce^{Ct}},
\end{equation}
where $C=C(\epsilon,\|u_0\|_{H^s(\mathbb{A}^d)})$>0.
Combining \eqref{eq7} with \eqref{3-6} shows that there exists a positive
constant $C=C(\epsilon, \|u_0\|_{H^s(\mathbb{A}^d)})$ such that
$$\|u(t,\cdot)\|_{H^3(\mathbb{A}^d)}\leq Ce^{Ce^{Ct}},\ \ \ \ \forall\ t\in[0,T),$$
which implies that \eqref{eq2-5} holds true when $1\leq d\leq 3.$ Therefore, we complete the proof of Lemma \ref{L3-4}.
\end{proof}

\begin{proof}[Proof of Theorem \ref{thm:main} \ref{thm:main:a}]
Suppose that the maximal existence time $T<+\infty$, then Lemma \ref{L3-4} implies that $\|u(t,\cdot)\|_{L^\infty(\mathbb{A}^d)}+\|\nabla u(t,\cdot)\|_{L^\infty(\mathbb{A}^d)}<+\infty$, $\forall\ t\in [0,T)$, which contradicts to Lemma \ref{L3-3}. Therefore, we complete the proof of Theorem \ref{thm:main} \ref{thm:main:a}.
\end{proof}

\section{Global attractor and asymptotic stability}
\newtheorem{remark4}{Remark}[section]
\newtheorem{theorem4}{Theorem}[section]
\newtheorem{lemma4}{Lemma}[section]
\newtheorem{corollary4}{Corollary}[section]
\newtheorem{proposition4}{Proposition}[section]

In this section, we will prove Theorem \ref{thm:main} \ref{thm:main:b} by deriving the continuous dissipative semigroup, the global attractor and global-in-time algebraic asymptotic stability. For this, let us firstly recall the definition of global attractor as follows. 

 \begin{definition1}\cite{MR953967}
 	Let $\{S(t)\}_{t\geq 0}$ be a semigroup of continuous operator on a Banach space $X.$ We say $A\subseteq X$ is a global attractor for the semigroup $\{S(t)\}_{t\geq0}$ in $X,$ if $A$ is compact, $S(t)A=A$, $\forall\ t>0$, and for any bounded set $B\subseteq X$, we have $\textup{dist}_X(S(t)B,A)\rightarrow0$, as $t\rightarrow+\infty$, where $\textup{dist}_X(Y,Z)\triangleq\sup\limits_{y\in Y}\inf\limits_{z\in Z}\|y-z\|_X$.
 \end{definition1}
 
For $k\in\mathbb{N}$, define $X_k(\mathbb{T}^d)=\{f\in H^k(\mathbb{T}^d)|\int_{\mathbb{T}^d}f(x)dx=0\}$. Then we have the following result.

\begin{lemma4}\label{lemma4-1}
Let $u_0\in X_3(\mathbb{T}^d)$, and 
$u(t,x)$ be the corresponding global solution to system \eqref{eq1}. Then $u(t,\cdot)\in X_3(\mathbb{T}^d), \ \forall \ t\geq0.$
\end{lemma4}

\begin{proof}
In view of Theorem \ref{thm:main} \ref{thm:main:a}, we have  $u(t,\cdot)\in H^3(\mathbb T^d).$ So, it suffices to show $\int_{\mathbb T^d}u(t,x)dx=0$. Indeed, since $\int_{\mathbb T^d}u_0(x)dx=0$, it follows that $\int_{\mathbb T^d}m_0(x)dx=0.$
By recalling the  system \eqref{eq1} in component, one has 
\begin{align*}
\partial_tm_i+u_j\partial_jm_i+\partial_ju_jm_i+m_j\partial_iu_j=\epsilon\Delta m_i.
\end{align*}
Then, integrating over $\mathbb T^d$ on both sides, we deduce
\begin{align*}
\frac{d}{dt}\int_{\mathbb T^d}m_i(t,x)dx&=-\int_{\mathbb T^d}(u_j\partial_jm_i+m_i\partial_ju_j+m_j\partial_iu_j)dx+\epsilon\int_{\mathbb T^d}\Delta m_idx\\
&=-\int_{\mathbb T^d}\left(\partial_j(u_jm_i)+(u_j-\Delta u_j)\partial_iu_j\right)dx=0,
\end{align*}
which implies that $\int_{\mathbb T^d}m(t,x)dx=\int_{\mathbb T^d}m_0(x)dx=0.$
By the relation $u=m+\Delta u$, we deduce that $\int_{\mathbb T^d}u(t,x)dx=0$. Therefore, for all $t\geq 0$, we have $u(t,\cdot)\in X_3(\mathbb{T}^d)$.
\end{proof}

\begin{remark4}
	If we denote $S(t)u_0(x)\triangleq u(t,x)$, then $S(t)$ is a semigroup of the solution operator for system \eqref{eq1}. Moreover, in view of Lemma \ref{L3-1}, we get that $S(t)$ is uniformly bounded in $H^1$ with $\|S(t)u_0\|_{H^1}\leq\|u_0\|_{H^1}$. 
\end{remark4}

\begin{remark4}\label{remark4-1}
 By Lemma \ref{lemma4-1}, for any $u_0\in X_3(\mathbb{T}^d)$, thanks to Poincar\'{e}'s inequality, we have  $\|u(t,\cdot)\|_{L^2}\leq C(d)\|\nabla u(t,\cdot)\|_{L^2}\leq C(d)\|\nabla^2u(t,\cdot)\|_{L^2}\leq C(d)\|\nabla^3u(t,\cdot)\|_{L^2}$, $\forall\ t\geq0.$
 Moreover, the Sobolev embedding theorem implies that
 \begin{align*} \|u(t,\cdot)\|_{L^\infty}\leq C(d)\|\nabla^2u(t,\cdot)\|_{L^2}\leq  C(d)\|\nabla^3u(t,\cdot)\|_{L^2},\quad \forall\ t\geq0.
\end{align*}
\end{remark4}

To prove the main theorem, we need another lemma as follows.

\begin{lemma4}\label{lemma4-2}\cite{MR953967}(The Uniform Gronwall Lemma)
 Let $g, h, y$ be three positive locally integrable functions on $(0, +\infty)$ such that $y'$ is locally integrable on $(0, +\infty)$, and 
 \[y'(t)\leq g(t)y(t)+h(t), \ \ \  \ \forall\ t\geq 0.  \]
If there holds
\[ \int_0^{t}g(s)ds\leq a_1,\ \ \ \int_0^{t}h(s)ds\leq a_2,\ \text{and}  \ \int_0^{t}y(s)ds\leq a_3,\ \  \ \forall\ t\geq 0,  \]
where $a_1,a_2,a_3$ are positive constants, then we have
\[ y(t)\leq\left(\frac{a_3}{t}+a_2\right)e^{a_1}, \ \forall\ t>0. \]
\end{lemma4}

\begin{lemma4}\label{lemma4-3}
	For every $f\in H^1(\mathbb T^d)$  with  $1\leq d\leq3$, we have
	\begin{align}\label{eq8}
		\|f\|_{L^4}\leq C(d)(\|f\|_{L^2}+\|f\|^{\frac{1}{4}}_{L^2}\|\|\nabla f\|^{\frac{3}{4}}_{L^2}).
	\end{align}
	Moreover, if $\int_{\mathbb T^d}f(x)dx=0$, then
	\begin{align}\label{eq9}
		\|f\|_{L^4}\leq C(d)\|\nabla f\|_{L^2}.
	\end{align}
\end{lemma4}
\begin{proof}
	By means of the interpolation inequality and the Sobolev embedding, we infer 
	\begin{equation*}
		\|f\|_{L^4}\leq \|f\|_{L^2}^{\frac{1}{4}}\|f\|_{L^6}^{\frac{3}{4}}\leq \|f\|_{L^2}^{\frac{1}{4}}\|f\|_{H^1}^{\frac{3}{4}},
	\end{equation*}
	which yields \eqref{eq8}. Thanks to the Poincar\'{e} inequality, \eqref{eq9} follows from \eqref{eq8}.
\end{proof}

So, in order to complete the proof of Theorem \ref{thm:main}, we only need to prove the following result.

\begin{theorem4}\label{the4-1}

Let  $u_{0}\in X_3(\mathbb{T}^d)$ and $S(t)u_0=u(t)$. Then the semigroup of the solution operator $\{S(t)\}_{t\geq0} $ for the system \eqref{eq1} has a global attractor $\{0\}$  in  $X_3(\mathbb{T}^d)$. Moreover, there exists a constant $C=C(d,\epsilon,\alpha,\|u_0\|_{H^3(\mathbb{T}^d)})>0$ such that 
\[ \|S(t)u_0\|_{H^3(\mathbb{T}^d)}\leq Ct^{-\frac{1}{2}},\quad\forall\ t>0. \]
\end{theorem4}

\begin{proof}
By virtue of Remark \ref{remark4-1}, we only need to estimate $\|\nabla m\|_{L^2}$. For this, we denote $y(t)\triangleq\int_{\mathbb T^d}|\nabla m|^2(t,x)dx$. If then follows from \eqref{3-7}, Lemma \ref{lemma4-1} and  \eqref{eq8}-\eqref{eq9} that
\begin{align}\label{eq12}
&\frac{1}{2}\frac{d}{dt}y(t)+\epsilon\int_{\mathbb{T}^d}|D^2 m|^2dx\nonumber\\
&\leq C\|\Delta m\|_{L^2}(\| u\|_{L^4}\|\nabla m\|_{L^4}+\|m\|_{L^4}\|\nabla u\|_{L^4})\nonumber\\
&\leq C\|D^2m\|_{L^2}(\|\nabla u\|_{L^2}\|\nabla m\|_{L^2}+\|\nabla m\|_{L^2}\|D^2 u\|_{L^2}+\|\nabla u\|_{L^2}\|\nabla m\|^{\frac{1}{4}}_{L^2}\|D^2 m\|^{\frac{3}{4}}_{L^2})\nonumber\\
&\triangleq I_1+I_2+I_3.
\end{align}
Thanks to the Young inequality \eqref{eq6}, we deduce
\begin{align}\label{eq10}
 I_1&=C\|D^2m\|_{L^2}\|\nabla u\|_{L^2}\|\nabla m\|_{L^2}\leq \frac{\epsilon}{4}\|D^2m\|_{L^2}^2+\frac{C^2\|\nabla u\|^2_{L^2}\|\nabla m\|^2_{L^2}}{\epsilon},
\end{align}
\begin{align}
 I_2&=C\|D^2m\|_{L^2}\|\nabla m\|_{L^2}\|\Delta u\|_{L^2}\leq \frac{\epsilon}{4}\|D^2m\|_{L^2}^2+\frac{C^2\|\nabla m\|^2_{L^2}\|\Delta u\|^2_{L^2}}{\epsilon},
\end{align}
 and 
\begin{align}\label{eq11}
 I_3&=C\|D^2m\|_{L^2}\|\nabla u\|_{L^2}\|\nabla m\|^{\frac{1}{4}}_{L^2}\|D^2m\|^{\frac{3}{4}}_{L^2}\nonumber\\
&\leq \|D^2m\|_{L^2}\left(\frac{\epsilon}{4}\|D^2m\|_{L^2}+\frac{C^{4}\|\nabla u\|^4_{L^2}\|\nabla m\|_{L^2}}{4(\frac{\epsilon}{3})^3}\right)\nonumber\\
&\leq\nonumber\frac{\epsilon}{4}\|D^2m\|_{L^2}^2+\frac{C^{4}\|D^2m\|_{L^2}\|\nabla u\|^4_{L^2}\|\nabla m\|_{L^2}}{4(\frac{\epsilon}{3})^3}\nonumber\\
&\leq\frac{\epsilon}{2}\|D^2m\|_{L^2}^2+\frac{1}{\epsilon}\left(\frac{C^{4}\|\nabla u\|^4_{L^2}\|\nabla m\|_{L^2}}{4(\frac{\epsilon}{3})^3}\right)^2.
\end{align}
Combining \eqref{eq12} with \eqref{eq10}-\eqref{eq11}, one yields
\begin{align}\label{4-10}
\frac{1}{2}\frac{d}{dt}y(t)\leq&\frac{C^2\|\nabla u\|^2_{L^2}\|\nabla m\|^2_{L^2}}{\epsilon}+\frac{C^2\|\nabla m\|^2_{L^2}\|\Delta u\|^2_{L^2}}{\epsilon}+\frac{1}{\epsilon}\left(\frac{C^4\|\nabla u\|^4_{L^2}\|\nabla m\|_{L^2}}{4(\frac{\epsilon}{3})^3}\right)^2\nonumber\\
=&\left(\frac{C^2\|\nabla u\|^2_{L^2}}{\epsilon}+\frac{C^2\|\Delta u\|^2_{L^2}}{\epsilon}+\frac{1}{\epsilon}\left(\frac{C^4\|\nabla u\|^4_{L^2}}{4(\frac{\epsilon}{3})^3}\right)^2\right)y(t)\nonumber\\
\triangleq&g(t)y(t),
\end{align}
where $g(t)\triangleq\frac{C^2\|\nabla u\|^2_{L^2}}{\epsilon}+\frac{C^2\|\Delta u\|^2_{L^2}}{\epsilon}+\frac{1}{\epsilon}\left(\frac{C^4\|\nabla u\|^4_{L^2}}{4(\frac{\epsilon}{3})^3}\right)^2$. Then by virtue of Lemma \ref{L3-1}, we have
\begin{equation*}
	\int^t_0\|\nabla u\|_{L^2}^8ds\leq \bigg(\sup_{t\in[0,T)}\|\nabla u\|^6_{L^2}\bigg)\int_0^t\|\nabla u\|_{L^2}^2ds\leq \bigg(\sup_{t\in[0,T)}\| u\|^6_{H^1}\bigg)\int_0^t\|\nabla u\|_{L^2}^2ds\leq \frac{1}{2\epsilon}\| u_0\|_{H^1}^8.
\end{equation*}
Thus, one gets
\begin{align}\label{4-11}
\int_0^tg(s)ds\leq\frac{C^2}{2\epsilon^2}\|u_0\|^2_{H^1}+\frac{1}{2\epsilon^2}\left(\frac{C^{4}\|u_0\|^4_{H^1}}{4(\frac{\epsilon}{3})^3}\right)^2\triangleq \frac{a_1}{2}.
\end{align}

On the other hand, thanks to \eqref{3-4}, Lemma \ref{lemma4-3} and the Young inequality \eqref{eq6} again, similar to \eqref{eq10} and \eqref{eq11}, we infer 
\begin{align*}
&\quad\frac{1}{2}\frac{d}{dt}\int_{\mathbb T^d}|m|^2dx+\epsilon\int_{\mathbb T^d}|\nabla m|^2dx\\
&\leq\nonumber3\|u\|_{L^4}\|m\|_{L^4}\|\nabla m\|_{L^2}\\
&\leq\nonumber C\|\nabla u\|_{L^2}(\|m\|_{L^2}+\|m\|^{\frac{1}{4}}_{L^2}\|\nabla m\|^{\frac{3}{4}}_{L^2})\|\nabla m\|_{L^2}\\
&\leq\nonumber\frac{\epsilon}{4}\|\nabla m\|^2_{L^2}+\frac{C^2\|\nabla u\|^2_{L^2}\|m\|^2_{L^2}}{\epsilon}+\frac{\epsilon}{2}\|\nabla m\|_{L^2}^2+\frac{1}{\epsilon}\left(\frac{C^4\|\nabla u\|^4_{L^2}\|m\|_{L^2}}{4(\frac{\epsilon}{3})^3}\right)^2,
\end{align*}
or what is the same,
\begin{align}\label{4-12}
\frac{d}{dt}\int_{\mathbb T^d}|m|^2dx+\frac{\epsilon}{2}\int_{\mathbb T^d}|\nabla m|^2dx&\leq\left(\frac{2C^2\|\nabla u\|^2_{L^2}}{\epsilon}+\frac{2}{\epsilon}\left(\frac{C^4\|\nabla u\|^4_{L^2}}{4(\frac{\epsilon}{3})^3}\right)^2\right)\|m\|^2_{L^2}\nonumber\\
&\triangleq f(t)\|m\|^2_{L^2},
\end{align}
where $f(t)\triangleq\frac{2C^2\|\nabla u\|^2_{L^2}}{\epsilon}+\frac{2}{\epsilon}\left(\frac{C^4\|\nabla u\|^4_{L^2}}{4(\frac{\epsilon}{3})^3}\right)^2$. If we denote $h(t)\triangleq\int_0^tf(s)ds,$ then $h$ is increasing. Thanks to Lemma \ref{L3-1} again, it follows that 
\begin{equation}\label{eq13}
	h(t)\leq \frac{C^2\|u_0\|^2_{H^1}}{\epsilon^2}+\frac{1}{\epsilon^2}\left(\frac{C^4\|u_0\|^4_{H^1}}{4(\frac{\epsilon}{3})^3}\right)^2\triangleq M_1, \ \forall\ t\geq0.
\end{equation}
Thus \eqref{4-12} implies that
\begin{align*}
\frac{d}{dt}(\|m\|^2_{L^2}e^{-h(t)})+\frac{\epsilon}{2}\|\nabla m\|^2_{L^2}e^{-h(t)}\leq0,
\end{align*}
which turns out that
\begin{align*}
\|m\|^2_{L^2}e^{-h(t)}+\frac{\epsilon}{2}\int_0^t\|\nabla m(s)\|^2_{L^2}e^{-h(s)}ds\leq\|m_0\|^2_{L^2}.
\end{align*}
Then we get
\begin{align*}
e^{-h(t)}\frac{\epsilon}{2}\int_0^t\|\nabla m(s)\|^2_{L^2}ds\leq\frac{\epsilon}{2}\int_0^t\|\nabla m(s)\|^2_{L^2}e^{-h(s)}ds\leq\|m_0\|^2_{L^2},
\end{align*}
which along with \eqref{eq13} yields that
\begin{align}\label{4-13}
\int_0^ty(s)ds\leq \frac{2}{\epsilon}\|m_0\|^2_{L^2}e^{h(t)}\leq\frac{2}{\epsilon}\|m_0\|^2_{L^2}e^{M_1}\triangleq a_3.
\end{align}
Since \eqref{4-10}, \eqref{4-11} and \eqref{4-13}, it follows from Lemma \ref{lemma4-2} that
\begin{align}\label{eq14}
y(t)\leq\frac{a_3}{t}e^{a_1}, \ \forall t>0.
\end{align}
In view of Remark \ref{remark4-1}, we get $\|u(t)\|^2_{H^3}\leq C(d)y(t)$.
Then \eqref{eq14} gives $$\|S(t)u_0\|^2_{H^3}\leq C(d)\frac{a_3}{t}e^{a_1}, \ \forall\ t>0.$$
So, $\lim\limits_{t\rightarrow+\infty}\|S(t)u_0\|_{H^3}=0,$ for all $u_0\in X_3(\mathbb{T}^d).$ That is, $\{0\}$ is the global attractor for $S(t)$ in $X_3(\mathbb{T}^d)$. Therefore, we complete the proof of Theorem \ref{the4-1}.
\end{proof}

\smallskip
\noindent{\bf Acknowledgments} 
The work of Guan is partially supported by National Natural Science Foundation of China under grants 12026216 and 12171493.  The work of Yan is partially supported by National Natural Science Foundation of China under grant 11971188.

	%\vspace{1\baselineskip} 

%	\noindent\textbf{Data Availability Statement.} This article has no associated data.

%	\vspace{1\baselineskip} 

%	\noindent\textbf{Declaration of generative AI use.} The authors report generative AI was not used in their research or preparation of this manuscript.

%	\vspace{1\baselineskip} 

%	\noindent\textbf{Conflict of interest.} The authors declare that there are no conflicts of interest.

\bibliographystyle{abbrv}
\bibliography{ref}

\end{document}